\magnification=\magstep1   
\input amstex
\UseAMSsymbols
\input pictex 
\vsize=23truecm
\NoBlackBoxes
\parindent=18pt
  
   \font\rmk=cmr8


\def\op{\text{\rm op}}

\def\Hom{\operatorname{Hom}}

\def\Ext{\operatorname{Ext}}

\def\rad{\operatorname{rad}}
\def\add{\operatorname{add}}

\def\Cok{\operatorname{Cok}}
\def\soc{\operatorname{soc}}

  \def\ss{\ssize }
\def\arr#1#2{\arrow <1.5mm> [0.25,0.75] from #1 to #2}

\def\s{\hfill \square} 
 
\def\injdim{\operatorname{inj.dim.}}
\def\projdim{\operatorname{proj.dim.}}
\def\Gdim{\operatorname{Gdim}}

\vglue.5truecm
\centerline{\bf Simple reflexive modules over finite-dimensional algebras}
                     \bigskip\medskip
\centerline{Claus Michael Ringel}
                \bigskip\medskip

\noindent {\narrower Abstract:  \rmk Let $\ss A$ be a finite-dimensional algebra.
If $\ss A$ is self-injective, then all modules are reflexive. Marczinzik recently has asked
whether $\ss A$ has to be self-injective in case all the simple modules are reflexive.
Here, we exhibit an 8-dimensional algebra which is not self-injective, but such that
all simple modules are reflexive (actually, for this example, the simple modules are
the only non-projective indecomposable modules which are reflexive).
In addition, we present some
properties of simple reflexive modules in general. Marczinzik had motivated his question
by providing large classes $\ss \Cal A$ of algebras such that any algebra in $\ss \Cal A$
which is not self-injective has simple modules which are not reflexive. However, as it turns
out, most of these classes have the property that any algebra in $\ss \Cal A$
which is not self-injective has simple modules which are not even torsionless.

	\medskip
\noindent 
\rmk Key words. Reflexive module, torsionless module. Self-injective algebra. 
	\medskip
\par}
	\bigskip
{\bf 1. Introduction.} 
	\medskip
Let $A$ be an artin algebra. 
For simplicity, we usually will assume that $A$ is a finite-dimensional $k$-algebra
where $k$ is an algebraically closed field.

The modules to be considered are usually 
left $A$-modules of finite length. Given a module $M$, let $M^* = \Hom(M,A)$ 
be its $A$-dual, and $\phi_M\:M \to M^{**}$ the canonical map from $M$ to $M^{**}.$
A module $M$ is said to be {\it torsionless} provided it is isomorphic to a submodule of a
projective module, or, equivalently, provided $\phi_M$ is injective.
And $M$ is {\it reflexive} provided $\phi_M$ is bijective. 
	
	\medskip
If $A$ is self-injective, then all modules are reflexive. Marczinzik recently has asked
whether $A$ has to be self-injective in case all the simple modules are reflexive.
The aim of this note is answer this question in the negative: we are going to
exhibit an 8-dimensional wild algebra whose simple modules are the only
indecomposable modules which are reflexive and not projective, see section 2.

In section 4 we present some
properties of simple reflexive modules in general. Since we
are interested in simple modules $S$ with $S^{**}$ isomorphic to $S$, we also
will look at the corresponding right modules $S^*$. Note that if
$P$ is a projective cover of $S$, then $S^*$ is isomorphic to a
submodule of the right module $P^*$. 
If $e$ is a primitive idempotent of $A$ with $P = Ae$, then $S^*$ can be identified
with the right ideal $e(\soc{}_AA) \subseteq eA$. In section 3, we collect some 
well-known characterizations of self-injective algebras for reference in sections 4 and 5.

Marczinzik had motivated his question
by providing large classes $\Cal A$ of algebras such that any algebra in $\Cal A$
which is not self-injective has simple modules which are not reflexive. 
We recall these results in section 5. As it turns
out, most of these classes have the property that any algebra in $\Cal A$
which is not self-injective has simple modules which are not even torsionless.
	\bigskip\bigskip
{\bf 2. Some examples.}
	\medskip
{\bf 2.1. First example.}
We are going to present {\it an 8-dimensional wild algebra $A$ which is not self-injective,
with 2 simple modules, both being reflexive} (actually, the simple modules turn out to be the only indecomposable modules which are reflexive and not projective).
	\medskip
Let $A$ be the path algebra of the quiver with two vertices $1$ and $2$, four arrows
$c\:1 \to 1,$ $b\:1 \to 2$ and $x,y,\:2\to 1$
$$
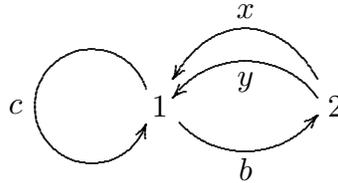
{\beginpicture
    \setcoordinatesystem units <1cm,1cm>
\multiput{} at 0 1  2 -0.5 /
\put{$1$} at 0.05 0
\put{$2$} at 2 0
\circulararc 320 degrees from -.1 0.2 center at -0.7 0
\ellipticalarc axes ratio 1:1 100 degrees from 1.8 0.1 center at 1 -.5
\ellipticalarc axes ratio 1:1 -100 degrees from 1.8 -.1 center at 1 .5
\ellipticalarc axes ratio 1:1.5 100 degrees from 1.8 0.3 center at 1 -.7
\put{$b$} at 1 -.7
\put{$c$} at -1.55 0
\put{$x$} at 1 1.05
\put{$y$} at 1 .27
\arr{-.13 -.27}{-.1 -.2}
\arr{1.75 -.15}{1.8 -.1}
\arr{0.25 .15}{.2 .1}
\arr{0.25 .38}{.2 .3}
\endpicture}
$$
and the relations
$$
  cx,\  by,\ 
  c^2,\ bc,\ 
  xb,\ yb.
$$
Thus, $A$ is a monomial algebra with basis 
$$
 e_1,\  e_2,\ x,\ y,\  b,\ c,\ bx,\ cy;
$$
the socle of ${}_AA$ has the basis $c,\ b,\ bx,\ cy.$

The two indecomposable projective modules can be visualized as follows:
$$
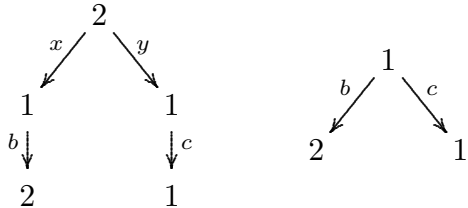
{\beginpicture
    \setcoordinatesystem units <.8cm,1cm>
\put{\beginpicture
\put{$2$} at 0 0
\put{$1$} at -1 -1
\put{$1$} at 1 -1
\put{$2$} at -1 -2
\put{$1$} at 1 -2
\arr{-0.2 -0.2}{-0.8 -.8}
\arr{0.2 -.2}{0.8 -.8}
\arr{-1 -1.3}{-1 -1.7}
\arr{1 -1.3}{1 -1.7}
\put{$\ssize x$\strut} at -0.6 -.3
\put{$\ssize y$\strut} at 0.6 -.3
\put{$\ssize b$\strut} at -1.2 -1.4
\put{$\ssize c$\strut} at 1.2 -1.4
\endpicture} at 0 0 
\put{\beginpicture
\put{$1$} at 0 0
\put{$2$} at -1 -1
\put{$1$} at 1 -1
\arr{-0.2 -0.2}{-0.8 -.8}
\arr{0.2 -.2}{0.8 -.8}
\put{$\ssize b$\strut} at -0.6 -.3
\put{$\ssize c$\strut} at 0.6 -.3
\endpicture} at 4 0 

\endpicture}
$$
	\medskip
The following modules will be of interest:
$$
 B = P(1)/S(1),\qquad C = P(1)/S(2).
$$
Note that $B, C$ are the indecomposable modules of length 2 with top $S(1)$. 

	\medskip 
{\bf (1)} {\it There are precisely six torsionless 
indecomposable modules, namely the projective modules $P(1), P(2)$, 
the simple modules $S(1), S(2)$ 
and the modules $B$ and $C.$}
	\medskip 
Proof. We have $J = \rad {}_AA = S(1)\oplus S(2) \oplus C \oplus B$. The annihilator of
$J$ is generated by $x,\ y$. A module $M$ is cogenerated by ${}_AJ$ if and only if 
$M$ is annihilated by $x,\ y$, thus the indecomposable modules cogenerated by ${}_AJ$
are $S(1), S(2), B, C, P(1).$ $\s$
	\medskip
{\bf (2)} {\it The $\mho$-quiver of $A$ consists of singletons and two components
of type $\Bbb A_3$, namely}
$$
{\beginpicture
    \setcoordinatesystem units <1.5cm,1cm>
\put{\beginpicture
\put{$S(1)$} at -.1 0
\put{$B$} at 1 0
\put{$\mho B$} at 2 0
\setdashes <1mm>
\arr{0.7 0}{0.3 0}
\arr{1.7 0}{1.3 0}
\endpicture} at 0 0 
\put{\beginpicture
\put{$S(2)$} at -.1 0
\put{$C$} at 1 0
\put{$\mho C$} at 2 0
\setdashes <1mm>
\arr{0.7 0}{0.3 0}
\arr{1.7 0}{1.3 0}
\endpicture} at 3.5 0 
\endpicture}
$$
	\medskip
Proof. We 
recall that the endpoints of the arrows of the $\mho$-quiver are just the 
indecomposable non-projective modules which are torsionless, see [RZ] 1.5. 
It is easy to see that the canonical inclusions $S(1) \to P(1),
S(2) \to P(1)$ are minimal left $\add{}_AA$-approximations, thus $\mho S(1) = B$
and $\mho S(2) = C.$ 
$\s$
	\medskip
The structure of the modules $\mho B$ and $\mho C$ is not really important, but since
some readers may be curious, we will provide some information below, see (6). Let us
mention already here: whereas $B$ and $C$ both are 2-dimensional (thus they look
rather similar), the module $\mho C$ is 3-dimensional, whereas $\mho B$ is 8-dimensional.
	\medskip
{\bf (3)} It follows that {\it $S(1)$ and $S(2)$ are the only indecomposable non-projective
modules which are reflexive,} see [RZ] 1.5.
$\s$
	\bigskip
It seems to be worthwhile to display also the opposite algebra $A^{\op}$ and to look at the
right modules $S(1)^*$ and $S(2)^*$. 
	\medskip
{\bf (4)} {\it The right modules $S(1)^*$ and $S(2)^*$ have length $2$.}
	\medskip
Proof. The socle of ${}_AA$ shows that 
$\dim S(i)^* = \dim\Hom(S(i),{}_AA) = 2$ for $i=1,2$.
$\s$
	\medskip
Here is the structure of the indecomposable projective right $A$-modules:
$$
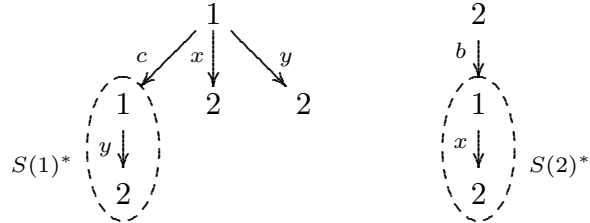
{\beginpicture
    \setcoordinatesystem units <1cm,1cm>
\put{\beginpicture
\put{$1$} at 0 0
\put{$1$} at -1 -1
\put{$2$} at 1 -1
\put{$2$} at 0 -1
\put{$2$} at -1 -2
\arr{-0.2 -0.2}{-0.8 -.8}
\arr{0.2 -.2}{0.8 -.8}
\arr{0 -.2}{0 -.8}
\arr{-1 -1.3}{-1 -1.7}
\put{$\ssize c$\strut} at -0.8 -.45
\put{$\ssize x$\strut} at -.18 -.45
\put{$\ssize y$\strut} at 0.8 -.45
\put{$\ssize y$\strut} at -1.2 -1.45
\setdashes <1mm>
\ellipticalarc axes ratio 1:2 360 degrees from -1 -.7 center at -1 -1.5
\put{$\ss S(1)^*$} at -1.9 -1.7
\endpicture} at 0 0 
\put{\beginpicture
\put{$2$} at 0 0
\put{$1$} at 0 -1
\put{$2$} at 0 -2
\arr{0 -.3}{0 -.7}
\arr{0 -1.3}{0 -1.7}
\put{$\ssize b$\strut} at -.2 -.4
\put{$\ssize x$\strut} at -.2 -1.4
\setdashes <1mm>
\ellipticalarc axes ratio 1:2 360 degrees from 0 -.7 center at 0 -1.5
\put{$\ss S(2)^*$} at .9 -1.7

\endpicture} at 4 0 
\endpicture}
$$
The radical $J_A$ (considered as a right $A$-module) is the direct sum of $S(1)^*$,
$S(2)^*$ and two copies of the simple right module with index $2$. 
We have encircled the submodules $S(1)^*$ and $S(2)^*$ by dashed lines.  
The right module $S(1)^*$ is the unique indecomposable submodule of $e_1A$ 
of length $2$. The right module $S(2)^*$ is the unique indecomposable submodule of $e_2A$ 
of length $2$. And we see:
	\medskip
{\bf (5)} {\it The right modules $S(1)^*$ and $S(2)^*$ are orthogonal bricks.}
	\bigskip
Let us rerun to the $\mho$-quiver of $A$. As we have mentioned, the structure of the modules
$\mho B$ and $\mho C$ is not really important, but it may be of interest to have them
available.
	\medskip
{\bf (6)} {\it The modules $\mho B$ and $\mho C$ can be displayed as follows:}
$$
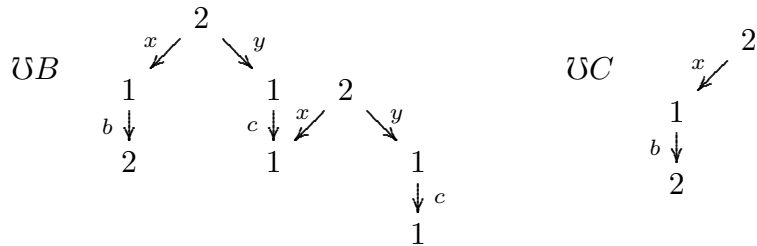
{\beginpicture
    \setcoordinatesystem units <.8cm,.8cm>
\put{\beginpicture
\put{$\mho B$} at -1.3 1.3 
\multiput{$1$} at 0 1  2 1  2 0  4 0  4 -1 /
\multiput{$2$} at 1 2  0 0  3 1 /
\arr{0 0.7}{0 0.3}
\arr{0.7 1.7}{0.3 1.3}
\arr{1.3 1.7}{1.7 1.3}
\arr{2 0.7}{2 0.3}
\arr{2.7 0.7}{2.3 0.3}
\arr{3.3 0.7}{3.7 0.3}
\arr{4 -.3}{4 -.7}
\multiput{$\ss x$\strut} at 0.3 1.7  2.4 0.7 /
\multiput{$\ss y$\strut} at 1.8 1.7  3.7 0.7 /
\multiput{$\ss b$} at -.3 0.5 /
\multiput{$\ss c$} at 1.7 0,5  4.3 -.5 /
\endpicture} at 0 0
\put{\beginpicture
\put{$\mho C$} at -1.2 1.6 
\multiput{$1$} at 0 1  /
\multiput{$2$} at 1 2  0 0   /
\arr{0 0.7}{0 0.3}
\arr{0.7 1.7}{0.3 1.3}
\multiput{$\ss x$\strut} at 0.3 1.7  /
\multiput{$\ss b$} at -.3 0.5 /
\endpicture} at 6 0.2
\endpicture}
$$

We have $\mho C = P(2)/Ax$ (note that $Ax$ is the only submodule of $P(2)$ isomorphic
to $C$), whereas $\mho B$ is an amalgamation of $P(2)$ and $P(2)/S(2).$ 
(As a first guess for $\mho B$, one may be inclined to 
consider $P(2)/Ax$, since $Ax \simeq B$,
but the diagram on the left shows that $\Ext^1(P(2)/Ax,P(2)) \neq 0$, whereas for any
module $M$, one has $\Ext^1(\mho M,{}_AA) = 0.$)
	\bigskip
{\bf 2.2.} {\bf A generalization.} {\it 
For any $n \ge 2$, there is a connected algebra $A = A(n)$
of dimension $2n+4$ with $n$ simple modules, all being reflexive, but $A$ is not
self-injective.}

For $n\ge 3$, let $A(n)$ be the path algebra of the quiver with vertices 
$1$, $2$, \dots $n$,
and arrows $x,y,\:n\to 1$, $c\:1 \to 1,$ $b\:1 \to 2$, and finally 
$a_i\:i \to i\!+\!1$ for $2\le i \le n\!-\!1$, thus:

$$
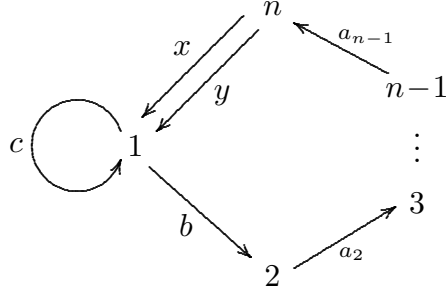
{\beginpicture
    \setcoordinatesystem units <.8cm,.8cm>
\multiput{} at 0 1  2 -0.5 /
\put{$1$} at 0.1 0
\put{$2$} at 2 -1.8
\put{$n$} at 2 1.9
\put{$n\!-\!1$} at 4 .8
\put{$3$} at 4 -.8
\circulararc 320 degrees from -.1 0.2 center at -0.7 0
\put{$b$} at .8 -1.1
\put{$c$} at -1.55 0
\put{$x$} at .75 1.35
\put{$y$} at 1.3 .7

\put{$\ss a_2$} at 3.1 -1.5
\put{$\ss a_{n-1}$} at 3.3 1.5

\arr{-.13 -.27}{-.1 -.2}

\arr{1.6 1.8}{0.2 0.4}
\arr{1.8 1.6}{0.4 0.2}

\arr{3.6 1}{2.3 1.7}
\put{$\vdots$} at 4 0.1
\arr{0.3 -0.3}{1.7 -1.55}
\arr{2.3 -1.7}{3.7 -0.85}
\endpicture}
$$
and the relations
$$
  cx,\  by,\ c^2,\ bc,\ a_2b,\ a_{i+1}a_i,\ xa_{n-1},\ ya_{n-1}.
$$
for $2\le i \le n\!-\!2$. 
Thus, $A$ is a monomial algebra with basis 
$$
 e_1,\dots, e_n,\ x,\ y,\  b,\ c,\ a_2,\dots, a_{n-1},\ bx,\ cy;
$$
and the indecomposable projective modules can be visualized as follows:
$$
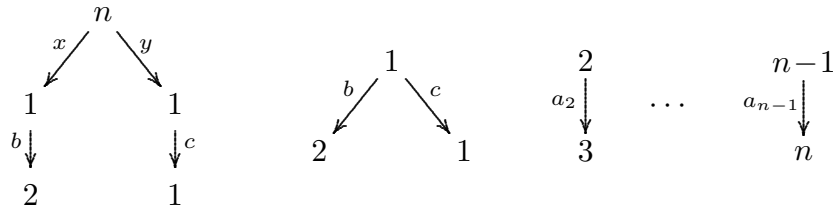
{\beginpicture
    \setcoordinatesystem units <.8cm,1cm>
\put{\beginpicture
\put{$n$} at 0 0
\put{$1$} at -1 -1
\put{$1$} at 1 -1
\put{$2$} at -1 -2
\put{$1$} at 1 -2
\arr{-0.2 -0.2}{-0.8 -.8}
\arr{0.2 -.2}{0.8 -.8}
\arr{-1 -1.3}{-1 -1.7}
\arr{1 -1.3}{1 -1.7}
\put{$\ssize x$\strut} at -0.6 -.3
\put{$\ssize y$\strut} at 0.6 -.3
\put{$\ssize b$\strut} at -1.2 -1.4
\put{$\ssize c$\strut} at 1.2 -1.4
\endpicture} at 0 0
\put{\beginpicture
\put{$1$} at 0 0
\put{$2$} at -1 -1
\put{$1$} at 1 -1
\arr{-0.2 -0.2}{-0.8 -.8}
\arr{0.2 -.2}{0.8 -.8}
\put{$\ssize b$\strut} at -0.6 -.3
\put{$\ssize c$\strut} at 0.6 -.3
\endpicture} at 4 0 
\put{\beginpicture
\put{$2$} at 0 0
\put{$3$} at 0 -1
\arr{0 -0.2}{0 -.8}
\put{$\ssize a_2$\strut} at -0.3 -.4
\endpicture} at 6.5 0 
\put{\beginpicture
\put{$n\!-\!1$} at 0 0
\put{$n$} at 0 -1
\arr{0 -0.2}{0 -.8}
\put{$\ssize a_{n-1}$\strut} at -0.45 -.4
\endpicture} at 9.5 0 
 
\put{$\cdots$} at 7.85 0 
\endpicture}
$$
	\medskip
Let 
$$
  B = P(1)/S(1),\quad C = P(1)/S(2), 
$$
then $B$ and $C$ are the indecomposable modules of length 2 with top $S(1)$.
	\medskip 
{\bf (1)} {\it There are precisely $2n+2$ torsionless 
indecomposable modules, namely the projective modules $P(i)$, 
the simple modules $S(i)$, with $1\le i \le n$ as well as $B$ and $C.$}
	\medskip 
{\bf (2)} {\it The $\mho$-quiver of $A$ consists of singletons and two components
of type $\Bbb A_3$ and $\Bbb A_{n+1}$ respectively, namely}
$$
{\beginpicture
    \setcoordinatesystem units <1.7cm,.8cm>
\put{\beginpicture
\put{$S(1)$} at -.05 0
\put{$B$} at 1 0
\put{$\mho B$} at 2 0
\setdashes <1mm>
\arr{0.7 0}{0.3 0}
\arr{1.7 0}{1.3 0}
\endpicture} at 1.1 0 
\put{\beginpicture
\put{$S(n)$} at -.2 0
\put{$\cdots$} at .9 0
\put{$S(2)$} at 1.95 0
\put{$C$} at 3 0
\put{$\mho C$} at 4 0
\setdashes <1mm>
\arr{0.6 0}{0.2 0}
\arr{1.6 0}{1.2 0}
\arr{2.7 0}{2.3 0}
\arr{3.7 0}{3.3 0}
\endpicture} at 0 -1 
\endpicture}
$$
	\medskip
{\bf (3)} {\it The simple modules are the only indecomposable non-projective modules 
which are reflexive.}
	\medskip
{\bf (4)} {\it The
right modules $S(1)^*$ and $S(2)^*$ have length $2$, 
the right modules $S(i)^*$ with $3\le i \le n$ are simple.}
	\medskip
Proof. The socle of ${}_AA$ shows that 
$\dim S(i)^* = \dim\Hom(S(i),{}_AA) = 2$ for $i=1,2$ and that 
$\dim S(i)^* = \dim\Hom(S(i),{}_AA) = 1$ for $3\le i \le n$.
$\s$
	\bigskip
The structure of the indecomposable projective right $A$-modules
is as follows:
$$
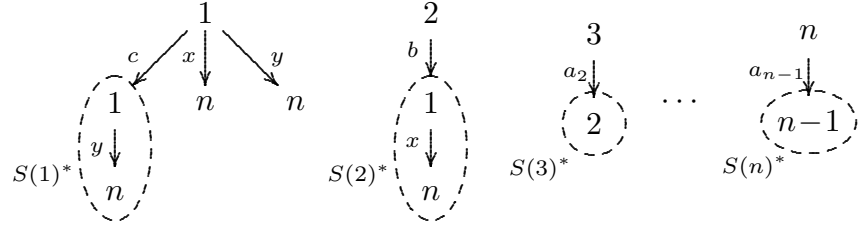
{\beginpicture
    \setcoordinatesystem units <1cm,1cm>
\put{\beginpicture
\put{$1$} at 0 0
\put{$1$} at -1 -1
\put{$n$} at 1 -1
\put{$n$} at 0 -1
\put{$n$} at -1 -2
\arr{-0.2 -0.2}{-0.8 -.8}
\arr{0.2 -.2}{0.8 -.8}
\arr{0 -.2}{0 -.8}
\arr{-1 -1.3}{-1 -1.7}
\put{$\ssize c$\strut} at -0.8 -.45
\put{$\ssize x$\strut} at -.18 -.45
\put{$\ssize y$\strut} at 0.8 -.45
\put{$\ssize y$\strut} at -1.2 -1.45
\setdashes <1mm>
\ellipticalarc axes ratio 1:2 360 degrees from -1 -.7 center at -1 -1.5
\put{$\ss S(1)^*$} at  -1.8 -1.8
\endpicture} at 0.6 0 
\put{\beginpicture
\put{$2$} at 0 0
\put{$1$} at 0 -1
\put{$n$} at 0 -2
\arr{0 -.3}{0 -.7}
\arr{0 -1.3}{0 -1.7}
\put{$\ssize b$\strut} at -.2 -.4
\put{$\ssize x$\strut} at -.2 -1.4
\setdashes <1mm>
\ellipticalarc axes ratio 1:2 360 degrees from 0 -.7 center at 0 -1.5
\put{$\ss S(2)^*$} at  -.8 -1.8
\endpicture} at 3.1 0 
\put{\beginpicture
\put{$3$} at 0 0
\put{$2$} at 0 -1
\arr{0 -.3}{0 -.7}
\put{$\ssize a_2$\strut} at -.2 -.4
\setdashes <1mm>
\ellipticalarc axes ratio 1:1 360 degrees from 0 -.65 center at 0 -1
\put{$\ss S(3)^*$} at  -.6 -1.5
\endpicture} at 5 0 
\put{\beginpicture
\put{$n$} at 0 0
\put{$n\!-\!1$} at 0 -1
\arr{0 -.3}{0 -.7}
\put{$\ssize a_{n-1}$\strut} at -.35 -.4
\setdashes <1mm>
\ellipticalarc axes ratio 1.5:1 360 degrees from 0 -.65 center at 0 -1
\put{$\ss S(n)^*$} at  -.6 -1.5
\endpicture} at 7.5 0 
\put{$\cdots$} at 6.4 0

\endpicture}
$$
The right module $S(1)^*$ is the unique indecomposable submodule of $e_1A$ 
of length $2$. The right module $S(2)^*$ is the unique indecomposable submodule of $e_2A$ 
of length $2$.
	\medskip
{\bf (5)} {\it The right modules $S(1)^*,$ \dots, $S(n)^*$ are pairwise orthogonal bricks.}
	\bigskip
{\bf 2.3.} Let us stress the following. Let $A$ be an algebra with radical $J$,
If $A$ is not self-injective, but all simple modules are reflexive, then $A$ cannot be 
local (see 5.3) and $J^2 \neq 0$ (see 5.4). 

	\bigskip\bigskip
{\bf 3. Characterizations of self-injective algebra using the socle.}
	\medskip
Let us collect some well-known characterizations of a self-injective algebra,
dealing with the socle of the algebra.
	\medskip
The algebra $A$ is said to be {\it left Kasch} provided and simple module is a submodule
of ${}_AA$. And $A$ is said to be {\it left QF-2} provided all indecomposable projective
modules are uniform (a uniform module of finite length is a module with simple socle).
	\medskip 
{\bf 3.1. Lemma.} {\it Let $A$ be a basic algebra with simple modules $S(1),\dots,S(n)$.
The following assertions are
equivalent:
\item{\rm (i)} $A$ is self-injective
\item{\rm (ii)} The modules $\soc P$ with $P$ indecomposable projective, are
    all the simple modules.
\item{\rm (iii)} The modules $S^*$ with $S$ simple, are all the simple right modules.
\item{\rm (iv)} $\soc {}_AA = \bigoplus_i S(i).$
\item{\rm (v)} $S(i)^*$ is a simple right module, for any $i$.
\item{\rm (v$'$)} $\dim S(i)^* = 1$, for any $i$.
\item{\rm (vi)} $S(i)^*$ is zero or a simple right module, for any $i$.
\item{\rm (vi$'$)} $\dim S(i)^* \le 1$, for any $i$.
\item{\rm (vii)} $\soc {}_AA$ is a submodule of $\bigoplus_i S(i).$
\item{\rm (viii)} $I({}_AA)$ is isomorphic to a submodule of $D(A_A).$
\item{\rm (ix)} $I({}_AA)$ is isomorphic to $D(A_A).$
\item{\rm (x)} $A$ is left Kasch and left QF-2.\par}
	\medskip
Some hints for the proof: One shows the implications (i) $\implies$ (ii) $\implies$
\dots $\implies$ (ix), and then (ix) $\implies$ (i). The implication
(viii) $\implies$ (ix) follows from  ${}_AA \subseteq I({}_AA)$, 
and $\dim A = \dim D(A).$ 

For the equivalence of (iv) and (x), one uses the fact that ${}_AA = \bigoplus P(j)$
with non-zero modules $P(j)$. Of course, (iv) implies that $A$ is left Kasch. 
Since $\soc P(j)$ is non-zero and $\soc {}_AA = \bigoplus
\soc P(j)$, it follows that $\soc P(j)$ is simple, for any $j$, thus $A$ is left QF-2. 
Conversely, if $A$ is left QF-2, then $\soc {}_AA = \bigoplus_j \soc P(j)$ is the direct sum
of $n$ simple modules. If $A$ is in addition left Kasch, then $\soc {}_AA$ is the direct
sum of the simple modules, each one occurring with multiplicity 1.
$\s$
	\bigskip
{\bf Remark.} Let us return to the examples exhibited in section 2. 
For these algebras, the right modules $S^*$ with $S$ simple are 
pairwise orthogonal bricks. 
	\bigskip
{\bf 3.2.} If all right modules $S^*$ with $S$ simple 
are simple  (condition (v)), then $A$ is self-injective,
thus all modules are reflexive, in particular all simple modules are reflexive.
The converse is not true, as we have seen in section 2.

Also, if $S$ is 
simple, then $S^*$ may be a simple right module, whereas $S$ is not
reflexive, see $S = S(2)$ for the quiver with relation
$$
{\beginpicture
    \setcoordinatesystem units <1cm,1cm>
\put{$1$} at  0 0
\put{$2$} at  1 0
\put{$\ 3\ .$} at  2 0
\arr{0.7 0}{0.3 0}
\arr{1.7 0}{1.3 0}
\setdots <.5mm>
\ellipticalarc axes ratio 1.5:1 130 degrees from 1.5  0.15 center at 1 0
\endpicture}
$$
Here, $S(2)$ is torsionless and $S(2)^{**} = P(3).$ 

Of course, if $S$ is simple, then $S^* \neq 0$ iff $\phi_S$ is injective thus iff
the kernel of $\phi_S$ is zero. 
But there is no obvious relationship between the dimension of
$S^*$ and the cokernel of $\phi_S.$
	\bigskip\bigskip
{\bf 4. Simple reflexive modules.}
	\medskip
Let us focus the attention to simple reflexive modules $S$ and their $A$-duals
$S^*$. If $A$ is commutative and $S$ is a simple module, then it is easy to see that
$S^*$ is semi-simple. As a consequence, for $A$ commutative and $S$ a simple reflexive
module, $S^*$ has to be again simple (this was stressed by Marczinzik, see [M], Proposition
3.11). In contrast, if $A$ is not commutative, then the $A$-dual $S^*$ of a simple module
may be a rather complicated right module, even if $S$ is reflexive. 
	\medskip
{\bf 4.1. Proposition.} {\it Let $S$ be a simple reflexive module. 
Then $S^*$ is a torsionless right module,  
and no proper non-zero factor module is torsionless. In particular, $S^*$
is a brick.}
	\medskip
Proof. Of course, $S^*$ is torsionless (see for example [RZ] 4.2). 
Assume there is a proper non-zero factor module
$N$ which is torsionless. Then there is an embedding $S^* \to (A_A)^m$ for some $m$
as well as a map $S^* \to (A_A)^{m'}$ for some $m'$ with image isomorphic to $N$.
This shows that $\dim S^{**} = \dim \Hom(S^*,A_A) \ge 2$, in contrast to the assumption 
$S^{**} = S$.

It follows that $S^*$ is a brick. Namely, if there is a non-zero non-invertible 
endomorphism, say with image $N$, then $N$ is a proper non-zero factor module which is torsionless.
$\s$ 
	\bigskip
{\bf 4.2. Proposition.} {\it Let $S, T$ be non-isomorphic 
simple reflexive modules. Then $S^*, T^*$
are orthogonal bricks.}
	\medskip
Proof. First, we show: 
{\it If $U,V$ are torsionless bricks and $\dim \Hom(U,A) = 1 = \dim\Hom(V,A),$
then $\Hom(U,V)\neq 0$ implies that $U$ and $V$ are isomorphic.} Namely, 
let $f\: U \to V$ be a non-zero map. Since $V$ is torsionless, there is a monomorphism
$V \to A^n$, thus a map $g\:V \to A$ with $gf\:U \to A$ non-zero. 
But $\dim\Hom(U,A) = 1$, say with $0\neq h\:U \to A$,
implies that the image of $gf$ is equal to $h(U) \simeq U$. Therefore, $f$ is a split monomorphism, thus an isomorphism. 

We apply this to $U = S^*, V = T^*$. Then $U,T$ are torsionless bricks. 
Note that $U,V$ are non-isomorphic: namely, $S^*$ and $T^*$ would be isomorphic, 
then $S = S^{**}$ and $T = T^{**}$ are isomorphic, a contradiction.
This shows that $\Hom(S^*,T^*) = 0 = \Hom(T^*,S^*).$ 
$\s$
	\bigskip
{\bf 4.3. Proposition.} {\it Assume that all simple modules are reflexive. If
$T$ is a simple module such that $S^*$ is a simple right module for all simple modules
$S \neq T$, then also $T^*$ is simple (and $A$ is self-injective).}
	\medskip
Proof. Let $S(1),\dots,S(n)$ be the simple modules and $T = S(1)$.
We assume that all the modules $S(i)$ are reflexive and that the
right modules
$S(2)^*,\dots, S(n)^*$ are simple. According to 3.2, the modules $S(1)^*,\dots,S(n)^*$
are pairwise orthogonal bricks. Thus, $S(2)^*,\dots,S(n)^*$ are pairwise different
simple right modules. Since there are $n$ simple right modules, let $N$ be the missing one.
Since $\Hom(S(1)^*,S(i)^*) = 0$ for $2\le i \le n$, all simple factor modules of $S(1)^*$
are isomorphic to $N$. Similarly, $\Hom(S(i)^*,S(1)^*) = 0$ for $2\le i \le n$ shows that
all simple submodules of $S(1)^*$ are isomorphic to $N$. Since $S(1)^*$ is a brick, it follows
that $S(1)^* = N,$ thus also $S(1)^*$ is a simple right module. 
It follows from 3.1 (iv) that $A$ self-injective.
$\s$
	\bigskip\bigskip
{\bf 5. Simple torsionless modules.}
	\medskip
Marczinzik has shown in [M] that there are several classes $\Cal A$ of artin algebras
such that any algebra $A$ in $\Cal A$ with all simple modules being reflexive, is
self-injective. We recall: $A$ is {\it Gorenstein} provided 
$\injdim {}_AA = \injdim A_A < \infty.$ Second, 
$A$ is  {\it QF-3} provided that $I({}_AA)$ is projective. Of course, as we have mentioned
already, $A$ is {\it left QF-2} 
provided any indecomposable projective module is uniform. 
The following entry collects nearly all relevant
assertions of [M].
	\medskip
{\bf 5.1. Marczinzik [M]:} {\it Assume that $A$ satisfies one of the following conditions.
	\smallskip
\item{\rm(a)} $\injdim {}_AA = \injdim A_A$ and $\projdim I({}_AA) < \infty$
  (for example, $A$ is Gorenstein). 

\item{\rm(b)} $\projdim I({}_AA) \le 1$ (for example, $A$ is QF-3).

\item{\rm(c)} $A$ is left QF-2.

\item{\rm(d)} Any simple right module is reflexive.
	\smallskip
\noindent
Then, if all simple (left) modules are reflexive, $A$ is self-injective.}
	\medskip 
Here are the references: (a) is [M] Theorem 3.6. (b) is [M] Proposition 3.8 
(and [M] Corollary 3.9).
(c) is [M] Proposition 3.10. (d) is [M] Proposition 3.15 (a). 
The paper [M] mentions two further classes $\Cal A$ such that any $A\in \Cal A$
with only reflexive simple modules is self-injective, namely the
artin algebras $A$ with $A^{\op}$ isomorphic to $A$ (it includes all commutative artin
algebras), see [M] Proposition 3.15 (b).
Now, if $A^{\op}$ is isomorphic to $A$ 
and any simple module is reflexive, then $A$ satisfies condition (d). 
Only the case of a local artin algebra is not covered by 5.1. 
	\medskip
As we will see, for all classes of algebras considered in 5.1,
one may replace the assumption that {\it the simple modules are reflexive} by the much weaker
assumption that {\it all simple modules are torsionless} and obtains the same conclusion:
	\medskip
{\bf 5.2. Proposition.} {\it Let $A$ be an artin algebra satisfying one of the conditions {\rm
(a), (b), (c), (d)} of 5.1. Then, if all simple (left) modules are torsionless, $A$
is self-injective.}
	\medskip
Proof (for (a), (b), (c), we follow quite closely [M], 
3.6, 3.7, 3.10, respectively).
We assume that $A$ is an artin algebra and that all simple modules are torsionless.
We may assume that $A$ is basic. 
	\medskip
(a) Assume that $\injdim {}_AA = \injdim A_A$ and $\projdim I({}_AA) < \infty$.
Since any simple module is a submodule
of ${}_AA$, any indecomposable injective module is a direct summand of $I({}_AA),$
thus $\projdim I({}_AA) < \infty$ implies that $\projdim D(A_A) < \infty.$
But $\projdim D(A_A) = \injdim A_A$, therefore $A$ is a Gorenstein algebra. 

Let $g = \injdim {}_AA$. We want to show that $g = 0.$
For the contrary, assume that $g \ge 1.$ For any module $M$, there is $0 \le t \le g$
such that $\Omega^tM$ is Gorenstein-projective and the smallest number $t$ is called
the Gorenstein-dimension $\Gdim M$ of $M$. It is well-known that always there exists 
a module of Gorenstein-dimension equal to $g$. If $S$ is simple and not projective,
then there is an indecomposable projective module $P$ and an embedding $u\:S \to P.$
Let $M = \Cok u,$ thus $\Omega M = S$.
Since $\Gdim M \ge g$, we see that $\Omega^{g-1}S = \Omega^g M$ is
Gorenstein-projective. This shows that all simple modules have Gorenstein-dimension at
most $g-1$. For any $t\ge 0$, the modules of Gorenstein-projective dimension 
at most $t$ are closed under extensions. Therefore, all modules have Gorenstein-dimension
at most $g-1$, a contradiction. $\s$
	\medskip
(b) We assume that $\projdim I({}_AA) \le 1$. According to [AR] Theorem 0.1,
this implies that the class of torsionless modules is closed under extensions.
Since all simple modules are torsionless, we see that all modules are torsionless.
In particular, any injective module is torsionless, thus a direct summand of a
projective module, thus projective. $\s$
	\medskip
(c) We assume that
$A$ is left QF-2, thus any indecomposable projective module is uniform. 
Let $S(1),\dots,S(n)$ be the simple modules.
Since all simple modules are torsionless, they occur in the socle of ${}_AA.$
Since ${}_AA$ is the direct sum of $n$ uniform modules, the socle of ${}_AA$
has length $n$, thus the socle of ${}_AA$ is equal to $\bigoplus_i S(i)$,
therefore $A$ is self-injective, see Lemma 3.1 (iv).
	\medskip
(d) We assume that any simple right module is torsionless. 
If $S$ is a simple module, then 4.1 asserts that $S^*$ is a right module
which has no proper non-zero factor module which is torsionless. Since $S^*$
is non-zero, there is a simple right module $N$ which is a factor module of $S^*$.
It follows that $N = S^*$, this shows that $S^*$ is simple. According to
Lemma 3.1 (v), $A$ is self-injective.
$\s$
	\medskip
We have mentioned above that 5.1 collects nearly all relevant results of [M], the only exception is the assertion that local artin algebras such that the simple 
module is reflexive, are self-injective (see [M] Proposition 3.14 (2)).
	\medskip 
{\bf 5.3 (Marzinzik).} {\it Let $A$ be a local algebra.
If all simple modules are reflexive, then $A$ is self-injective.}
	\medskip
Proof. Let us stress that also 5.3 is 
an immediate consequence of Proposition 5.2. Namely, if $A$ is a local artin algebra and its
simple module is reflexive, then the algebra $A^{\op}$ satisfies the condition (d):
the only simple right $A^{\op}$-module is reflexive and the only simple left $A^{\op}$-module
is of course torsionless. Thus, $A^{\op}$, and therefore $A$, is self-injective. 
	\bigskip
A further result of this kind should be mentioned.
	\medskip
{\bf 5.4. Proposition.} {\it Let $A$ be an algebra with radical square zero.
If all simple modules are reflexive, then $A$ is self-injective.}
	\medskip
Proof. We assume that $A$ is connected. 
Let $S$ be a simple module. If $S$ is injective, then $S$ is also
projective (since $S$ is torsionless), thus $A$ is a simple algebra and therefore
self-injective. 

Thus, we can assume that there are no simple injective modules (this means: the quiver of
$A$ has no source). Let $S$ be a simple module. Since $S$ is not injective, let $T$
be a simple module with $\Ext^1(T,S) \neq 0.$ Then $S$ is a submodule of $\rad P(T)$.
If $S$ is a proper submodule of $\rad P(T)$, then $\mho S$ has length at least 2,
(and is indecomposable and not projective). But this implies that $\mho S$ is not torsionless,
since torsionless modules are projective or simple. This contradiction shows that
$S = \rad P(T)$ and therefore $\dim\Ext^1(T,S) = 1$. Similarly, if there are non-isomorphic
simple modules $T,T'$ wish $\Ext^1(T,S) \neq 0 \neq \Ext^1(T',S)$, then $\mho S$
has length at least 2 (and is indecomposable and not projective...). Altogether, we see that
$A$ is a radical square zero Nakayama algebra whose quiver is a cycle, thus $A$
is self-injective.
$\s$
	\bigskip\bigskip
{\bf 6. References.}
	\medskip 
\item{[AR]} M. Auslander, I. Reiten. Syzygy modules for noetherian rings. J\. Algebra
   183 (1996), 167--185.
\item{[M]} R. Marczinzik. Simple reflexive modules over Artin algebras.
   Journal of Algebra and its Applications. (2019) 1950193.
\item{[RZ]}  C. M. Ringel, P. Zhang. 
   Algebra \& Number Theory. Vol. 14 (2020), No. 1, 1-36,  DOI: 10.2140/ant.2020.14.1
	\bigskip\bigskip
{\baselineskip=1pt
\rmk
C. M. Ringel\par
Fakult\"at f\"ur Mathematik, Universit\"at Bielefeld \par
POBox 100131, D-33501 Bielefeld, Germany  \par
ringel\@math.uni-bielefeld.de}

	\bigskip
\bye